\documentclass[12pt]{article}
\usepackage{bm}
\usepackage{epsf}
\usepackage{amssymb}
\usepackage{graphics}
\usepackage{amsthm}
\usepackage{amsfonts}
\usepackage{latexsym}
\usepackage{mathrsfs}
\usepackage{graphicx}
\usepackage{amsxtra}
\usepackage{mathtools}

\usepackage{anyfontsize}
\usepackage{t1enc}

\nonstopmode
\usepackage[T2A]{fontenc}
\usepackage{inputenc}
\usepackage[russian,english]{babel}
\usepackage{caption}
\inputencoding{cp1251}
\usepackage{morewrites}
\usepackage{etoolbox}

\usepackage[usenames,dvipsnames]{color}
\usepackage{hyperref}
 \hypersetup{pdfborder={0 0 0},
   colorlinks,
   urlcolor={blue},
   linkcolor={blue},
   citecolor={blue},
   breaklinks=true}
\usepackage{color}
\usepackage{xcolor}

\usepackage{doi}
\providecommand{\eprint}[1]{}
\renewcommand{\eprint}[1]{arXiv:\href{http://arxiv.org/abs/#1}{#1}}

\usepackage{titlesec}
\titleformat{\chapter}[display]
{\normalfont%
    \normalfont
    \bfseries}{\chaptertitlename\ \thechapter}{20pt}{%
    \Large 
    }

\usepackage[textwidth=165mm,textheight=235mm]{geometry}

\newcommand{\beqn}{\begin{eqnarray}}
\newcommand{\eeqn}{\end{eqnarray}}

\newcommand{\beq}{\begin{equation}}
\newcommand{\eeq}{\end{equation}}
\newcommand{\bs}{\bigskip}

\newcommand{\veps}{\varepsilon}
\newcommand{\HH}{\mathbb{H}}
\newcommand{\op}{\mathrm{op}}
\newcommand{\pa}{\partial}

\newtheorem{theorem}{Theorem}

\title{The true story of Quantum Ergodic Theorem}

\author{Alexander Shnirelman}

\date{} 

\begin{document}
{\large

\maketitle

Let $M$ be a smooth compact Riemannian manifold, and $\Delta$ the Laplace--Beltrami operator on $M$. Let $u_n$ and $\lambda_n$ be the eigenfunctions and eigenvalues of $\Delta$, i.e. $\Delta u_n+\lambda_n u_n=0$. Let $\Omega=\{(x,\xi)\in T^* M, |\xi|=1\}$ be the bundle of unit covectors, and $g_t: T^* M\to T^* M$ be the geodesic flow in $T^* M$. Then $\Omega$ is invariant under $\{g_t\}$, and $\omega=\displaystyle\frac{dx\wedge d\xi}{d |\xi|}$ is the $(2 n-1)$-form on $\Omega$ invariant under $\{g_t\}$ and thus defining the $g_t$-invariant (Liouville) measure $\mu$.

\begin{theorem}[\mbox{\cite{shnirelman1974ergodic,shnirelman1993,zelditch1987uniform,colin1985ergodicity}}]
Suppose the geodesic flow $\{g_t\}$ is ergodic with respect to the measure $\mu$. Then there exists a subsequence $u_{n_k}$ of eigenfunctions having density 1 such that for any smooth function $\varphi\in C^\infty(M)$, 
$$
\lim_{k\to\infty} \int_M \varphi(x) |u_{n_k}(x)|^2\,dx =\int_M \varphi(x)\,dx.
$$
\end{theorem}

\bigskip
\bigskip

Here I'll tell the discovery
story of this theorem. This story consists of a number of steps, some of them
being not so obvious.

\bs

{\bf 1.} My teacher Mark Iosifovich Vishik was one of the founders of the
Microlocal Analysis and Pseudodifferential Operators. He developed his
original version of it, which in some respects (PDO in bounded domains,
factorization) went far beyond the standard theory. His Mehmat seminar was
a crucible of the new concepts in this domain. No wonder that it was in
M.I.'s seminar where Egorov explained for the first time his celebrated
theorem, and I could learn it firsthand (though I appreciated its importance
much later).

\bs

{\bf 2.} When I was in the 3rd year of Mehmat studies, M.I. assigned me the first
research topic, namely the solvability (or Fredholmness) of singular integral
(i.e. 1-dimensional pseudodifferential of order zero) equations degenerating
at one point. To my surprise, I managed to find a satisfactory solution of this
problem (Fredholmness condition and appropriate functional spaces). In the
process I discovered for myself a (rather primitive) version of the wavelet
decomposition of a function, i.e. its microlocalization in the phase space.

\bs

{\bf 3.} In the 4th year student work I studied the difference equations in the bounded domain. Using the discrete version of factorization, I found a correct formulation of Boundary Value Problems in the convex domains.  

In my 5th year graduate thesis (which
became my first published work) I found new classes of Fredholm convolution equations with constant symbol in half-space.

My PhD. thesis was
devoted to topological methods in nonlinear problems of complex analysis,
and had nothing to do with eigenfunctions. It had no resonance during the
next 50 years, though some ideas may be interesting even today.

\bs

{\bf 4.} While in the graduate school, I bumped into the work
of Babich and Lazutkin
\cite{babich1967eigenfunctions}
and of Lazutkin \cite{lazutkin1968construction-r}
(translated in \cite{babich1968eigenfunctions,lazutkin1968construction}).
In these works, they constructed asymptotic solutions to the 
Helmholtz equation on a closed surface concentrated near a closed and stable geodesic
(Figure~\ref{qet-figure-1}\,(a))
and  in a bounded domain
concentrated near a closed billiard trajectory
(for example, the forth-and-back trajectory on Figure~\ref{qet-figure-1}\,(b)).
It should be noted that at that time the foreign
journals were almost unavailable for us. It was in
part because the subscription was quite limited, and
in part because of our poor English and other
languages except Russian. Practically we were
confined to the Russian journals, and a few foreign
articles given to us by our teachers. However, we were able to read all the
Russian journals on display in the Mehmat library, including obscure
ones. The works of Babich and Lazutkin were published in  such
unassuming places.
But they were true gems! The first consequence of my reading of those articles was
related to the famous work of Arnold, ``Modes and Quasimodes'' \cite{arnold1972modes}. 

\begin{figure}
\begin{center}
(a)
\includegraphics[width=0.25\textwidth]{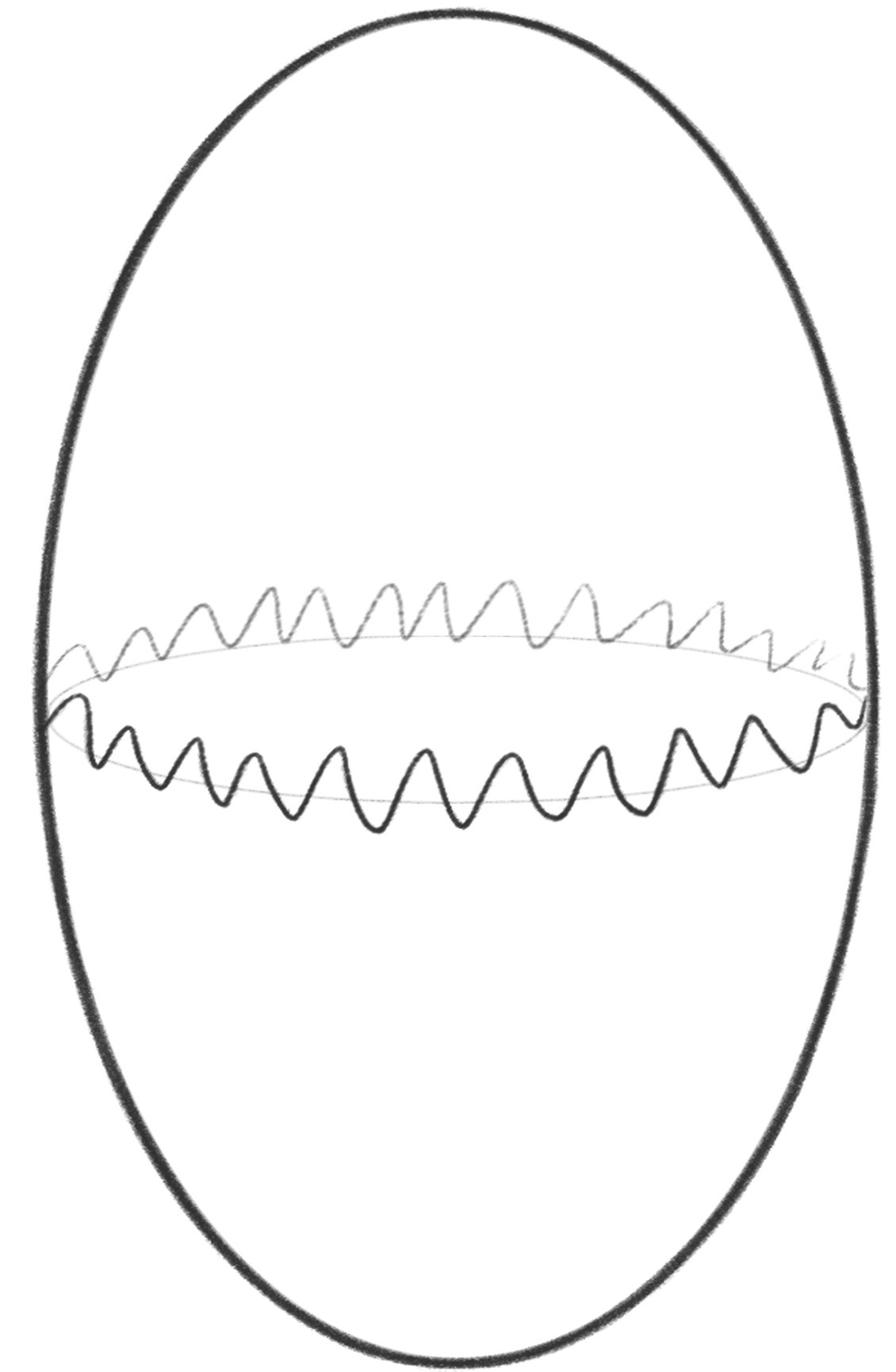}
\hskip 3.3cm
(b)
\includegraphics[width=0.45\textwidth]{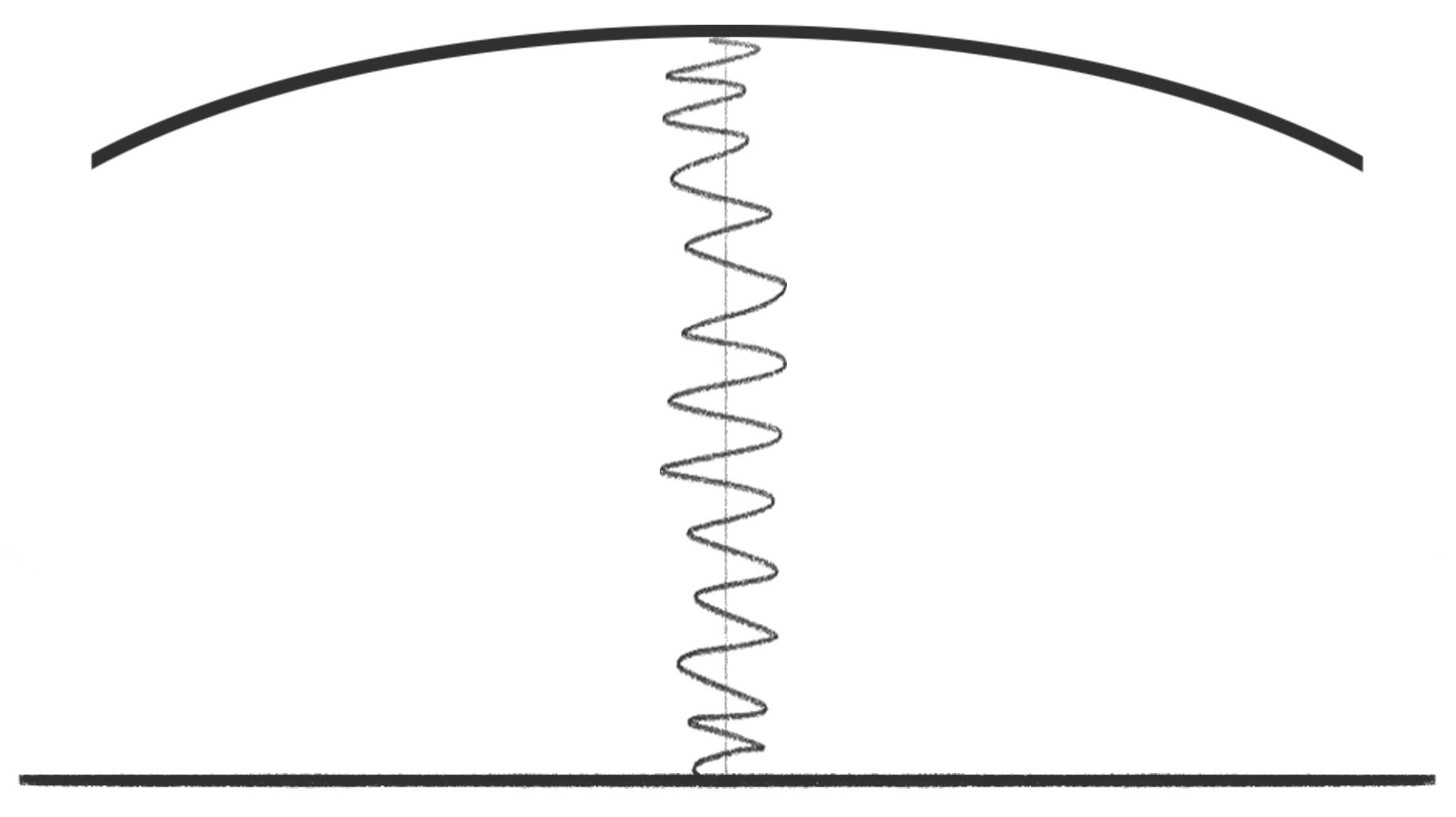}


\caption{
(a)
An eigenfunction on a closed surface concentrated near a closed (dynamically stable) geodesic.
(b)
An eigenfunction in a bounded domain concentrated near a bouncing trajectory.
}
\label{qet-figure-1}
\end{center}\end{figure}

Once Arnold told in his seminar about his new results on the spectrum of
a linear oscillating system with symmetries. He argued that if, for example,
the system has a symmetry of the 3rd order, then a part of the eigenvalues
have a stable multiplicity 2, and the remaining eigenvalues are simple. In
general, there are no eigenvalues of stable multiplicity 3: they split into
simple and double ones (according to the dimensions of irreducible
representations of the group of order 3).
 At this point I stood up and boldly
declared that this statement is wrong, and that
I have a counterexample. Arnold asked me to
produce it. So, I drew on the blackboard an
equilateral triangle with rounded corners, and
said that according to the results of  Lazutkin, such a membrane (with fixed
boundary) has ``laser'' eigenfunctions
concentrated near three altitudes of the
triangle
(Figure~\ref{qet-figure-2}).
They do not feel one another, and are
transformed one into another by a rotation by 120 degrees.
Hence the eigenvalue is stably triple.
Arnold scratched his head, and did not find what to answer. However, when I
met him in a week, he was shining with joy. ``I know what happens here!
These eigenfunctions are just approximate! In fact they overlap a little, and
so anyway there is a small (exponential) splitting into the simple and double
eigenvalues''. (In fact I've pretty much confused things. The works of
Babich and Lazutkin are rather difficult, and there are some arithmetic
conditions deep inside them making an obstacle to the continuous deformation of
eigenfunctions in the course of a change of a domain.)
 Afterwards Arnold described his theory in his famous article ``Modes
and Quasimodes'' \cite{arnold1972modes}. In this paper he describes this story in detail. Nevertheless,
this example is known as ``Arnold's example'', in spite of the explicit
attribution in his article.

\begin{figure}
\begin{center}
\includegraphics[width=0.4\textwidth]{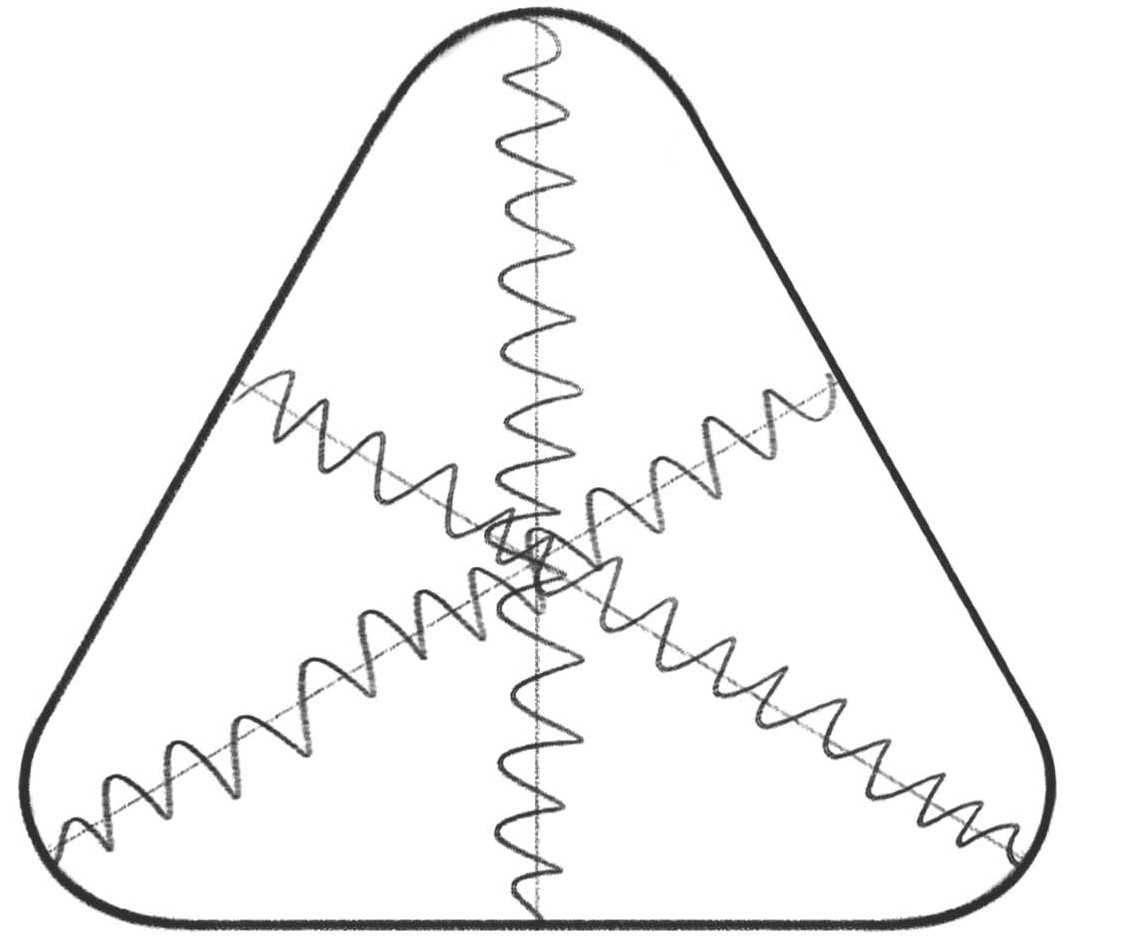}
\caption{
Laser-type eigenfunctions in a domain with the 3rd order symmetry;
an almost-counterexample to Arnold's conjecture.
}
\label{qet-figure-2}
\end{center}\end{figure}

\bs

{\bf 5.} Since my first reading of the work of Babich and Lazutkin
\cite{babich1968eigenfunctions,lazutkin1968construction},
I could not
think of anything else but the high-frequency eigenfunctions.
In \cite{babich1968eigenfunctions} and \cite{lazutkin1968construction},
the
construction hinges on the assumption that the billiard trajectory is
dynamically stable. And what about the unstable trajectories? I tried to
reproduce their  construction  (which, in the first approximation, reduced to
some ansatz resulting in the quantum harmonic oscillator equation in the
transverse direction to the trajectory); of course, I failed. The same scheme resulted in the Weber equation describing the scattering of waves on the potential barrier; it could not produce any discrete eigenvalues and corresponding asymptotic eigenfunctions concentrated near an unstable periodic trajectory.

\bs

{\bf 6.} Gradually I came to the following general formulation: ``How do the high-frequency eigenfunctions look like in general?'' I did not have a clear idea what
the word ``look'' exactly means. And there was a good reason for such
fuzziness. My teacher Mark Iosifovich Vishik, before his deep works on the
microlocal analysis, devoted several years to the joint work with Lazar
Aronovich Lyusternik on the asymptotic behavior of solutions to elliptic
equations with a small parameter at the higher order terms (what is called ``singular perturbations'') \cite{vishik1957regular}. Of course, I was aware of this activity. The problem looked quite similar to the high-frequency
eigenfunction problem which, too, can be regarded as a problem with a small
parameter (namely, inverse of the eigenvalue) at the higher order term
(Laplacian). The difference was the sign of the small parameter: it was
negative for the boundary layer situation, and positive for the eigenfunction
problem. Hence, we have a clear and simple picture of asymptotic behavior
in one case (a smooth core and thin boundary layers \cite{vishik1957regular}; later internal layers
were added to the picture after the work of Ventsel and Freidlin \cite{ventsel1970small}), and a hell
in the other. The (now classical) asymptotic methods used by Vishik and Lyusternik gave nothing for the study of high-frequency eigenfunctions. 

{\bf 7.} There was a well-developed asymptotic theory of eigenfunctions for the classically integrable systems, i.e. for Riemannian manifold $M$ such that the geodesic flow on the surface $\Omega\subset T^* M$ is completely integrable \cite{keller1958corrected,fedoryuk1976quasiclassical}. A typical example is the high-frequency eigenfunction in a disk with Dirichlet's condition on the boundary (Figure~\ref{qet-figure-3}). Later this theory was generalized by Lazutkin to the Riemannian manifolds whose geodesic flow is a small perturbation of the integrable one. The book \cite{lazutkin1993kam} contains exhaustive results in this direction. However, this theory gives no clue to the understanding of high-frequency eigenfunctions in the general, non-integrable case.

\begin{figure}
\begin{center}
\includegraphics[width=0.7\textwidth]{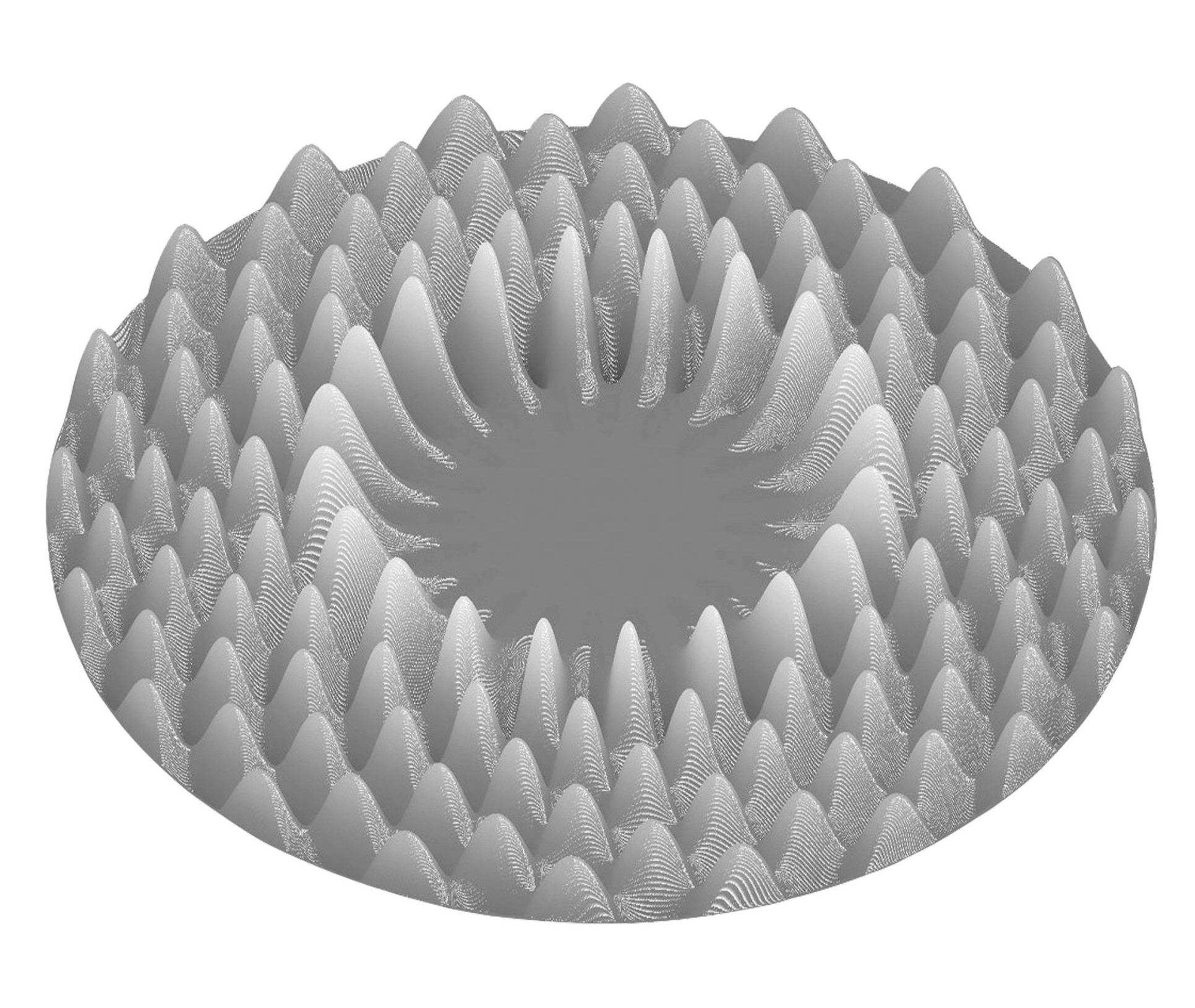}
\caption{
High-frequency eigenfunction
($J_n(j_{n,k}r)\sin(n\theta)$ with $n=20$ and $k=6$)
in a unit disk with Dirichlet's condition on the boundary.
%
}
\label{qet-figure-3}
\end{center}\end{figure}

\bs

\begin{figure}
\vskip -0.55cm
\begin{center}
\includegraphics[width=0.4\textwidth]{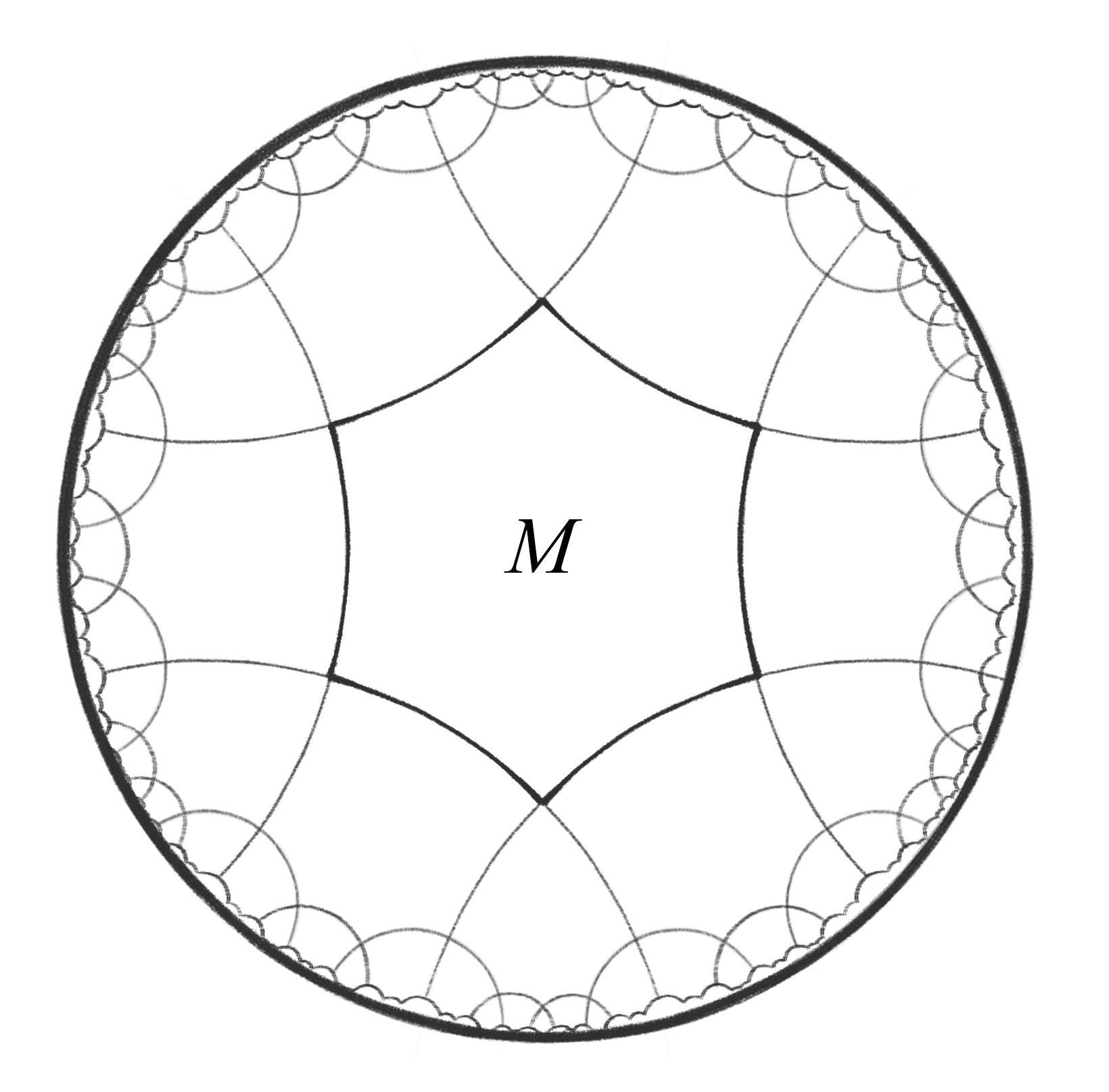}
\caption{
The surface $M$ as a fundamental domain of the discrete group $\Gamma$.
}
\label{qet-figure-4}
\end{center}\end{figure}

{\bf 8.} I decided to consider the case which is maximally remote from the
classically integrable cases, like the one studied by  Lazutkin. Of course, it was the
Laplacian on a compact surface of constant negative curvature. Why
constant? Because, I thought, it may give us some extra tools to study the
problem. In this I was right, and at the same time profoundly wrong; after all,
this restriction was superfluous, and hence
misleading.
So, I considered the surface which was a factor of
the action of a Fuchsian group $\Gamma$  of discrete
transformation of the hyperbolic plane   $\HH$:  $M=\HH/\Gamma$ (Figure~\ref{qet-figure-4}).
If $u_n$
is an eigenfunction on $M$, $\Delta u_n + \lambda_n u_n=0$,
then it can be lifted to a $\Gamma$-invariant (or automorphic) eigenfunction $v_n$
 on $\HH$: $v_n(x)=v_n(gx)$ ($x\in\HH, \  g\in \Gamma$). 

\bs

{\bf 9.} My first idea was to find a
representation of an automorphic
eigenfunction $v_n$
 similar to the
Poincare series for automorphic forms
and functions. Soon I realized that the
natural building blocks of such
representation are the ``horospheric'' 
eigenfunctions $w_n(x,y)$, i.e. such that their level curves are the horospheres
touching the absolute at one point $y$ (Figure~\ref{qet-figure-5}).
If $\rho$ is the distance
from a horosphere
to a fixed ``zero'' one, then $w_n(x,y)=e^{(-\frac{1}{2} + i k_n)\rho}$, and $\lambda_n=k_n^2+\frac{1}{2}$.

\bs

\begin{figure}
\begin{center}
\includegraphics[width=0.4\textwidth]{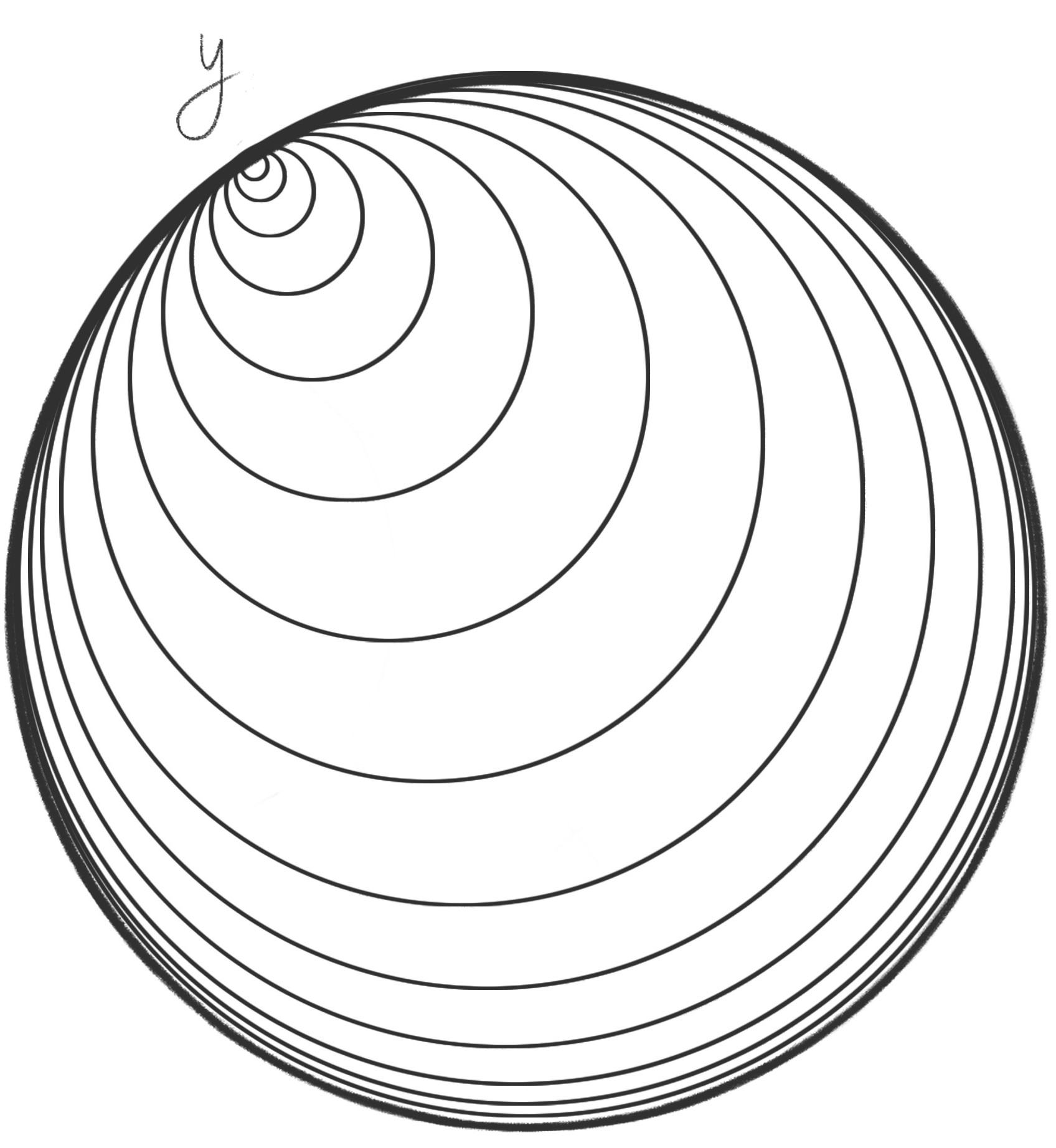}
\caption{
Level curves of a horospheric eigenfunction
$W_n(x,y)$ ($x\in\HH$, $y\in A$).
}
\label{qet-figure-5}
\end{center}\end{figure}

{\bf 10.} My first achievement was the proof that for every automorphic 
eigenfunction $v_n(x)$ there exists a distribution  $\Phi_n(y)$ of order not
exceeding $1$ on the absolute $A$ such that 

\beq
v_n(x)=\int_A \Phi_n(y) w_n(x,y)\,dy.
\eeq
Later on I learned that this statement is known under the name
``Helgason's Theorem'' \cite{helgason2018spherical},
and is the central fact of the harmonic analysis on locally symmetric spaces.

\bs

{\bf 11.} The next step was the functional (or homological) equation for the function
$\Phi_n(x,y)$. Namely, for any  $g\in\Gamma$,

\beq\label{qet-2}
\Phi_n(g(y))=|g'(y)|^{-1/2}e^{i k_n g'(y)}\Phi_n(y).
\eeq

Thus, we have an overdetermined system of functional equations on one
single function $\Phi_n(y)$. Fortunately, we do not need to care about its
solvability: we know a priori that it has solutions for the discrete sequence of
numbers  $k_n=\sqrt{\lambda_n-1/2}$. So, we can concentrate on the study of properties of the function $\Phi_n(y)$.

\bs

{\bf 12.} The action of the group $\Gamma$
 on the absolute $A$ is quite complicated; all the orbits are dense, and
there is no invariant measure. So, it is worthwhile to lift the action of $\Gamma$ onto
the cotangent space $T^* A$.  This action already preserves the Liouville
measure $dy\wedge d\eta$, but (a) for almost all points $(y,\eta)\in T^* A$ their orbits
are dense, and (b) the volume of the phase space $T^* A$ is infinite.

\bs

{\bf 13.} My goal was to prove that the functions $\Phi_n$ for all, or at least almost all $n$ look like typical realizations of the white noise. This property is best expressed in terms of the {\it Wigner measures} $\mu_n$ corresponding to the functions $\Phi_n$. For every function $\Phi_n$, we can define the {\it Wigner measure} $\mu_n$ in $T^* A$. This measure describes the distribution of the energy of the function $\Phi_n(y)$ in the phase space. It is defined in the following way. For any smooth function (symbol) $a(y,\eta)$ with compact support, let $\hat a$ be the corresponding (Weyl) pseudodifferential operator with the symbol $a$ \cite{hormander1994analysis3}. Now consider the bilinear expression $(\hat a \Phi_n, \Phi_n)$. This expression linearly depends on the symbol $a$, i.e. it has the form
\[
(\hat a\Phi_n,\Phi_n)=\int_{T^* A} U_n(y,\eta) a(y,\eta)\,dy\,d\eta.
\]
It turns out that the distribution $U_n(y,\eta)$ is asymptotically (as $n\to\infty$) a non-negative measure which we denote by $\mu_n$; this measure is called {\it Wigner measure}. 

\bs

{\bf 14.}
There is a physical device which produces the Wigner measure of a signal (i.e.
function of one variable; say, the time). It consists of a number of filters cutting a
narrow frequency band from the signal; each filter is characterized by the
median frequency. If we plot the squared filter outputs as a function of time
and the median frequency, we get a positive and highly oscillating function
(sometimes called the Husimi function). After some mollifying we get a
smooth positive density, which is the density of the Wigner measure. It is conjectured that our actual hearing has a similar mechanism, i.e. our brain converts the incoming sound into its Wigner measure which is further processed to extract the meaningful information.

 Some cases of the Wigner measure have been known since long ago: I mean the music scores. Consider a musical sound, say a song played by some instrument or sung by a human voice. The corresponding music score is a 2-dimensional domain $\bm{\mathit{Score}}$ endowed with two coordinates (time and frequency) with the notes which we regard as points. Let us put into each such point a mass proportional to the intensity of the corresponding note; then we get a measure on
the $\bm{\mathit{Score}}$ which is a coarse approximation
of the true Wigner measure $\mu$ in the time-frequency plane.

\bs

{\bf 15.}
The significance of the Wigner measure for our problem was based on the following observation. The quantity $I=\int_M f(x) |v_n(x)|^2\,dx$ can be transformed (after some manipulations) into $I=(B_f \Phi_n, \Phi_n)$ where $B_f$  is a pseudodifferential operator with a smooth and compactly supported symbol $b_f(y,\eta)$ depending on $f$ (I do not write this symbol explicitly, but its most important property is that $\int_{T^* A}b_f(y,\eta)\,dy\,d\eta=\int_M f(x)\,dx$). Thus,
\[
\int_M f(x) |v_n(x)|^2\,dx \sim (B_f\Phi_n,\Phi_n)\sim\int_{T^* A}b_f\,d\mu_n.
\]
So, in order to prove the asymtotic equidistribution of $v_n$, we have to prove that the measures $\mu_n$ are asymptotically equivalent to the Liouville measure $\Lambda=dy\wedge d\eta$ on $T^* A$.

\bs

{\bf 16.} Then I was able to prove that  the ``arithmetic mean'' of the measure $\mu_n$ is equal to the Liouville measure $\Lambda$. The exact meaning of this result is that $4\pi\tau \sum_n e^{-\lambda_n\tau}\mu_n \sim \Lambda$ as $\tau\to 0$.
It was done in the traditional way, with the use of the heat equation on $M$ (Carleman's method 
\cite{berline2004heat}). I have to confess that at that time I did not master the Tauberian theorem, and therefore did not make the next step, and only proved that
$\frac{1}{N}\sum_{n=1}^N \mu_n\sim \Lambda$ for $N\to\infty$.
This ignorance can be seen in the first publications of my results; later this gap was filled \cite{shnirelman1993}.

\bs

{\bf 17.} The relations (\ref{qet-2})
imply that for any $g\in\Gamma$, the measure $\mu_n$ is asymptotically (as $n\to\infty$) invariant under the transformation

\beq
(y,\eta) \mapsto F_g(y,\eta)=\big(g(y),(g'(y))^{-1}\eta+k_n g''(y)\big).
\eeq

\bs

{\bf 18.} Then I proved that the transformations $F_g$ defined above have the property of equidistribution (in the sense of Kazhdan \cite{kazhdan1965uniform}). To define this property, introduce the word distance $d(g,h)$ in the group $\Gamma$ as the minimal length $p$ of the word $i_1,\,\dots,\,i_p$ such that $g\cdot h^{-1}=\gamma_{i_1}^{\pm 1}\cdot\,\ldots\,\cdot \gamma_{i_p}^{\pm 1}$, where $\gamma_i$ are generators of the group $\Gamma$. Let $B_R=\{g\in \Gamma\, |\, d(g, e)\le R\}$. The family of transformations $\{F_g\}, g\in \Gamma$, possesses the equidistribution property if for any $z_0\in T^* A$ and any two bounded domains $U,\,V\subset T^* A$, 

\beq
\lim_{R\to\infty}\frac{|\{g\in B_R|F_g(z_0)\in U\}|}{|\{g\in B_R | F_g(z_0)\in V\}|}=\frac{|U|}{|V|}.
\eeq

\bs

{\bf 19.} The above properties of the functions $\Phi_n$ imply that for a sequence $n_k$ of density 1, the Wigner measures $\mu_{n_k}$ tend weakly to the Liouville measure
$\Lambda=dy\wedge d\eta$
in $T^* A$ (I've used the version of the Birkhoff ergodic theorem). And this implies the conclusion of Theorem 1, i.e. the asymptotic equidistribution of a subsequence $v_{n_k}$ of density 1. This completed the proof for the eigenfunctions of the Laplacian on the compact surface $M=\HH/\Gamma$ having constant negative curvature.

\bs

{\bf 20.} At this point I realized that it was possible to define the Wigner measures for the functions $v_k$ themselves as measures on $T^* M$, and to work with them, thus skipping a big part of the proof. The Wigner measures are defined in a manner similar to the above 1-d definition. Namely, for any smooth, compactly supported  function (symbol) $a(x,\xi)$  ($(x,\xi)\in T^*M$), let $\hat a = \op(a)$ be the pseudodifferential operator with the symbol $a$. Consider the bilinear expression $(\hat a v_n,v_n)$. Using the G\aa rding inequality
\cite{hormander1994analysis3},
we show that there exists a sequence of nonnegative measures $\mu_n$ on
$T^* M$ of mass 1 such that $(\hat a v_n, v_n)\sim\int_{T^* M} a\,d\mu_n$.

\bs

{\bf 21.}
Here I have to make some remarks on the G\aa rding inequality
\cite{hormander1994analysis3}.
Let $a(x,\xi)$ be a symbol of the H\"ormander class $S^0_{1,0}$, i.e. $|\pa^\alpha_x \pa^\beta_\xi a|\le C_{\alpha,\beta}(1+|\xi|)^{-|\alpha|}$; let $\hat a$ be the pseudodifferential operator with the symbol $a(x,\xi)$. Then the G\aa rding inequality (in fact, a pair of inequalities) says that if the symbol $a(x,\xi)\in S^0_{1,0}$ is real and nonnegative, then 
\begin{eqnarray}
{\rm Re} (\hat a u,u)&\ge& -C \Vert u\Vert^2_{H^{-1/2}},
\\[1ex]
|{\rm Im}(\hat a u,u)|&\le& C \Vert u\Vert^2_{H^{-1/2}},
\end{eqnarray}
where $H^{-1/2}$ is the Sobolev space.

Here the constant $C$ is common for all symbols belonging to a bounded set in the space $S^0_{1,0}$. Examples of symbols of class $S^0_{1,0}$ are (a) symbols homogeneous in $\xi$ of degree zero, and (b) symbols of the form $a(x,\xi)=\alpha(x,\veps \xi)$, where $\alpha\in C_0^\infty$ and $\alpha(x,\xi)=0$ for $\xi$ close to 0 (such symbols belong to a bounded set in $S^0_{1,0}$ uniformly for all $0<\veps<1$). So, when I refer to the G\aa rding inequality, I mean the symbols of type (b).

\bs

{\bf 22.} There exists a natural device which recovers the Wigner measure for
oscillating functions in any dimension. This device is our eye! If the eye is
put at some point in the space, then it will see some brightness distribution on
the ``sky'', and in every direction it will see some colour. Thus we have a
measure (the energy distribution) depending on the position, direction, and
frequency, i.e. exactly in the phase space. A closer analysis of the work of the eye (or any similar optical device) shows that, in fact, the visible picture seen by the eye on the ``sky'' is nothing but the Wigner measure. For example, suppose that the eigenfunction $u_n$ is ``quasiclassical'', i.e. if locally
\begin{equation}\label{qet-eigenfunction}
u_n(x)=\sum_{k=1}^K a_k(x) e^{i \lambda_n^{1/2}\varphi_k(x)}.
\end{equation}
Here the phase functions $\varphi_k(x)$ satisfy the equation $|\nabla \varphi_k|\equiv 1$, and the amplitudes $a_k(x)$ satisfy the transport equation $\nabla\varphi_k\cdot\nabla|a_k|^2+\Delta\varphi_k\cdot |a_k|^2=0$.  The points $(x,\nabla\varphi_k(x))\in T^* M$ form locally a Lagrangian manifold $\Lambda_k$; for different $k$, these manifolds are called \emph{Lagrangian sheets}.
The ``eye'' put at the point $x$ will ``see'' $K$
discrete stars on the dark ``sky'' in the directions $\nabla\varphi_k, \ \ k=1,\,\dots,\,K$. 

\bs

{\bf 23.}
Some pairs of Lagrangian sheets can merge along an $(n-1)$-dimensional manifold called {\it caustic}. These two sheets, together with the natural projection $\pi: T^* M \to M$, form a singularity called a \emph{fold} in the singularity theory.
If two Lagrangian sheets of the eigenfunction \eqref{qet-eigenfunction},
$a_k(x) e^{i \lambda_n^{1/2} \varphi_k(x)}$ and $a_m(x) e^{i \lambda_n^{1/2}\varphi_m(x)}$, merge together along a caustic,
and the ``eye''  is moving and its trajectory crosses the caustic, it will ``see'' two ``stars'' approaching one another, and then merging together, forming a single, much brighter ``star''. However, upon further motion of the ``eye'', it will see that the ``star'' almost immediately fades down and disappears. This exactly happens for the high-frequency eigenfunctions of the Laplacian inside the disk showed on
Figure~\ref{qet-figure-3}.
There is a central domain where the eigenfunction
is of order $e^{-C\sqrt{\lambda_n}}$
(the shadow domain). It is bounded by the concentric circle, the caustic. On this circle, the amplitude of the eigenfunction attains its maximum. In the annulus between the caustic and the boundary of the disk the eigenfunction is ``quasiclassical'', i.e. has the form
\beq
u_n(x)=\sum_{k=1}^2 b_k(x)\cos(\sqrt{\lambda_n}\varphi_k(x)).
\eeq
So, if the eye is in the ``quasiclassical'' domain, it sees two pairs of stars (one pair is opposite to the other) on the dark ``sky'' (here the ``sky'' is a circle);
see Figure~\ref{qet-figure-6}\,(a).
As the eye moves and approaches the caustics, these pairs of stars merge and form two stars at the opposite positions on the sky, whose brightness becomes much higher (Figure~\ref{qet-figure-6}\,(b)). But if the eye moves further into the shadow domain, the stars fade away, and the eye enters
the darkness (Figure~\ref{qet-figure-6}\,(c)).

\begin{figure}
\begin{center}
\includegraphics[width=0.6\textwidth]{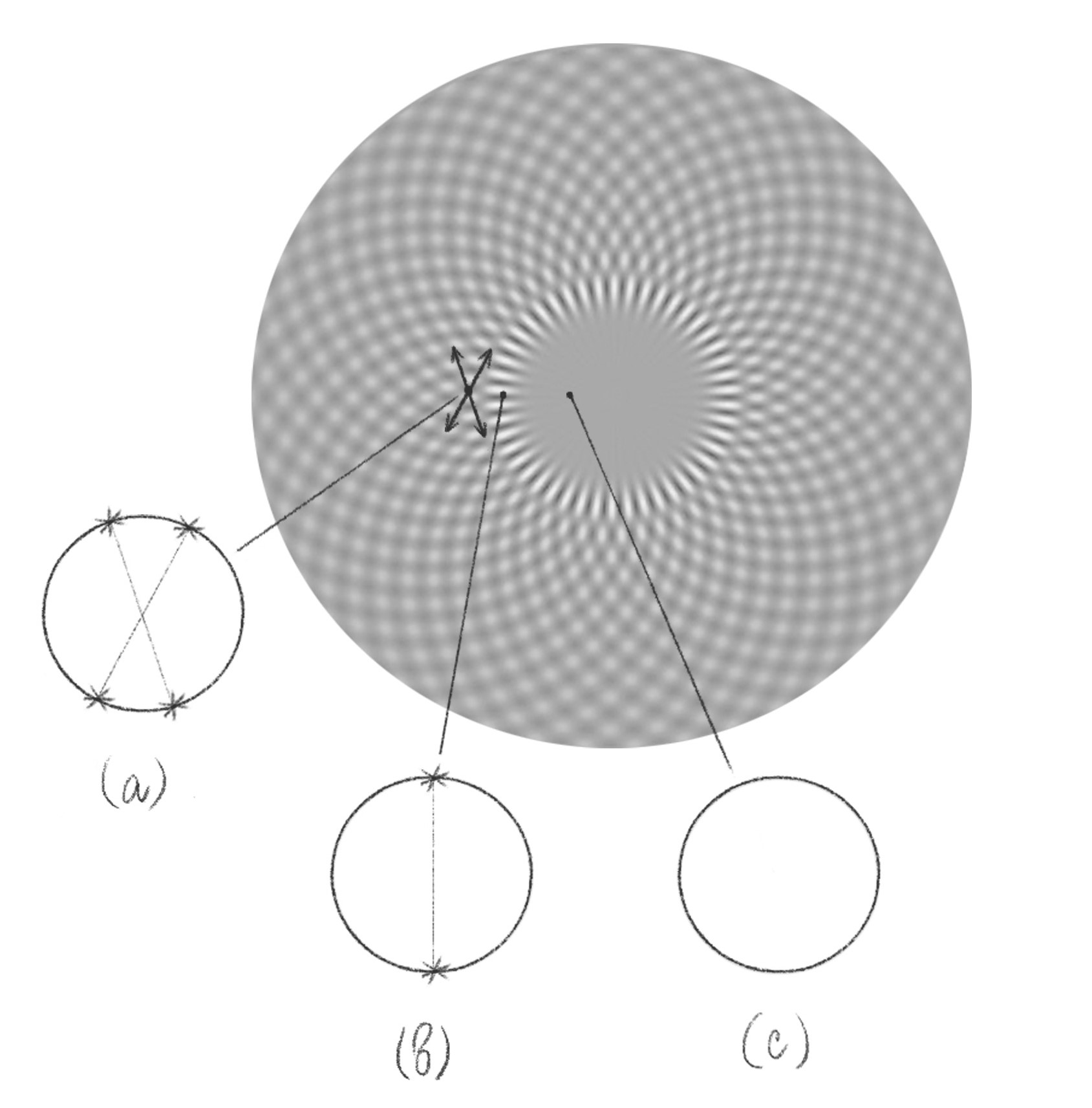}
\caption{
What the eye sees on the sky
from the quasiclassical domain (a),
from a point on the caustic (b),
and from the shadow domain (c).
The directions to the four starts on firmament (a)
are given by normals to the two wave fronts
of the wave function at a given point of the domain
(these directions are shown on the picture).
On the caustic, four directions degenerate into two.
}
\label{qet-figure-6}
\end{center}\end{figure}

\bs

{\bf 24.} Further, we prove that 

\begin{enumerate}
\item
Each measure $\mu_n$ in $T^* M$ is concentrated near the hypersurface $|\xi|=\sqrt{\lambda_n}$;
\item
$\sum_n \mu_n\sim dx\wedge d\xi$;
\item
The measures $\mu_n$ are asymptotically invariant under the geodesic flow (at this point I used the hyperbolic equation $u_t=i\sqrt{-\Delta}u$ and the  Egorov Theorem \cite{hormander1994analysis3,hormander1994analysis4});
\item
If the measures $\mu_n$ have the properties (1)--(3),
then there exists a subsequence $\mu_{n_k}$ of density 1 such that the measures $\mu_{n_k}$ are asymptotically equidistributed on the energy surfaces $|\xi|=\lambda_{n_k}^{1/2}$; the proof is based on the Birkhoff ergodic theorem.
\end{enumerate}

\bs

{\bf 25.} 
Returning to our optical interpretation, our result means that if
we consider the typical high-frequency eigenfunction $u_{n_k}$
(where $\{n_k\}$ is the subsequence of density 1 for which the equidistribution
of the Wigner measure holds), and put the eye at any point of the manifold $M$,
it will see the sky uniformly glowing; the intensity of this glow does not
depend on the eye's position and direction of its sight. Thus, the properties
of the eigenfunctions are quite opposite to the quasiclassical ones:
in the quasiclassical case, an eigenfunction is concentrated
on an $n$-dimensional Lagrangian manifold,
while in the ergodic case, an eigenfunction is uniformly distributed over a
finite energy submanifold of dimension $2n-1$
(this value can be interpreted as the dimension of the space of light rays in the cotangent bundle).

{\bf 26.} This proof of quantum ergodicity appears quite different
from the previous one,
sketched above in paragraphs {\bf 9}\,--\,{\bf 19}.
In place of the discrete group $\Gamma$ acting on the measures $\mu_n$ by (3) we have a  continuous one-parameter group of sympectic
transformations of the phase space (the phase flow). Thus, in the first proof, we have a discrete non-commutative
group of symplectic transformations of the space $T^* A$, the cotangent bundle of the absolute $A$ of the Lobachevsky plane, while in the second proof we have a continuous 1-parameter group of symplectic transformations of the space $T^* M$, the cotangent bundle of the surface $M$. The properties of the Wigner measures in these two proofs are quite different, too: in the first proof they are spread over the whole phase space, for the functions $\Phi_n(y)$ are singular even for small $n$, and in the second proof each measure $\mu_n$ is concentrated on the energy surface $|\xi|=\lambda_n^{1/2}$. These (and other) differences hint at the possibility that in these two proofs different structures were used, and this difference can result in
future interesting results. This opinion is confirmed by the excellent achievement of Rudnik and Sarnak \cite{rudnick1994behaviour} who proved that on the arithmetic surface {\it all} eigenfunctions are asymptotically equidistributed (this property is called Quantum Unique Ergodicity). These authors used an additional structure present in the arithmetic case, namely the existence of the Hecke correspondences, providing extra symmetries to the Wigner measures. 

\bs

{\bf 27.} Upon the proof of the Quantum Ergodicity Theorem it  immediately became clear that 

\begin{enumerate}
\item
The curvature of the surface $M$ can be non-constant (we need only ergodicity of the phase flow) \cite{zelditch1987uniform}; 
\item
The dimension of $M$ may be arbitrary \cite{colin1985ergodicity};
\item
We can consider any elliptic operator, not only Laplacian \cite{colin1985ergodicity};
\item
The manifold $M$ may have a boundary \cite{zelditch1996ergodicity}; 
\item
We can consider not only scalar, but also vector- and
bundle-valued functions and matrix operators on them \cite{gerard1991microlocal,gerard1993ergodic};
\item
There appear the first impressive results in the direction of Quantum Unique Ergodicty \cite{anantharaman2008entropy}.
\end{enumerate}

\bs

{\bf 28.} This was the Past. Now comes the Future.

\bigskip
\bigskip
\bigskip

\qquad
Alexander Shnirelman

\qquad
Concordia University, Montreal, Canada

\qquad
alexander.shnirelman@concordia.ca

}

\end{document}